\documentclass[12pt]{amsart}
\usepackage{amssymb,amsmath}
\makeatletter
\def\@cite#1#2{{\m@th\upshape\bfseries%
[{#1\if@tempswa{\m@th\upshape\mdseries, #2}\fi}]}}
\makeatother
\theoremstyle{plain}
\newtheorem{thm}{Theorem}[section]
\newtheorem{lem}[thm]{Lemma}
\newtheorem{cor}[thm]{Corollary}
\newtheorem{prop}[thm]{Proposition}
\theoremstyle{definition}
\newtheorem{rem}[thm]{Remark}
\newtheorem{defn}[thm]{Definition}
\newtheorem{note}[thm]{Note}
\newtheorem{eg}[thm]{Example}

\newtheorem{prob}[thm]{Problem}
\newcommand{\Prf}{\noindent\textbf{Proof.\ }}
\newcommand{\bx}{\strut\hfill$\blacksquare$\medbreak}
\newcommand{\upbx}{\vspace{-2.5\baselineskip}\newline\hbox{}\hfill$
\blacksquare$\newline\medbreak}

\newcommand{\ca}{\mathrm{C}^*}

\newcommand{\ol}{\overline}

\newcommand{\td}{\widetilde}

\newcommand{\sot}{\textsc{sot}}
\newcommand{\wot}{\textsc{wot}}

\newcommand{\bbB}{{\mathbb{B}}}
\newcommand{\bbC}{{\mathbb{C}}}
\newcommand{\bbD}{{\mathbb{D}}}
\newcommand{\bbE}{{\mathbb{E}}}
\newcommand{\bbF}{{\mathbb{F}}}

\newcommand{\bbR}{{\mathbb{R}}}

\newcommand{\bbU}{{\mathbb{U}}}

 \newcommand{\B}{{\mathcal{B}}}

\renewcommand{\H}{{\mathcal{H}}}

\renewcommand{\L}{{\mathcal{L}}}
 \newcommand{\M}{{\mathcal{M}}}

\renewcommand{\phi}{\varphi}
\newcommand{\upchi}{{\raise.35ex\hbox{$\chi$}}}
\newcommand{\fA}{{\mathfrak{A}}}

\newcommand{\fL}{{\mathfrak{L}}}

\newcommand{\fR}{{\mathfrak{R}}}

\newcommand{\qand}{\quad\text{and}\quad}
\newcommand{\qwhere}{\quad\text{where}\quad}

\newcommand{\qfor}{\quad\text{for}\quad}

\newcommand{\qif}{\quad\text{if}\quad}

\newcommand{\Alg}{\operatorname{Alg}}

\newcommand{\Lat}{\operatorname{Lat}}

\newcommand{\spn}{\operatorname{span}}

\newcommand{\sumin}{\sum_{i=1}^n}

\newcommand{\bofh}{\B(\H)}
\newcommand{\bofhn}{\B(\H_n)}
\newcommand{\rowa}{(A_1, \ldots, A_n)}
\newcommand{\rowt}{(T_1, \ldots, T_n)}
\newcommand{\rows}{(S_1, \ldots, S_n)}
\newcommand{\rowl}{(L_1, \ldots, L_n)}

\newcommand{\fnplus}{{\bbF}_n^+}

\newcommand{\basisvector}{\xi_w}
\newcommand{\vacuumvector}{\xi_e}
\newcommand{\lamiw}{\lambda_{i,w}}
\newcommand{\muiw}{\mu_{i,w}}

\newcommand{\flambda}{\fL_{\Lambda}}
\newcommand{\frightlambda}{\fR_{\Lambda}}
\newcommand{\nulam}{\nu_\lambda}
\newcommand{\wlam}{w(\lambda)}

\begin{document}

\title[weighted shifts on Fock space
]%
{Non-Selfadjoint operator algebras generated by weighted shifts on Fock space}
%
\author[D.W.Kribs]{David~W.~Kribs${}^1$}
\thanks{{\it key words and phrases. } Hilbert space, weighted shift, left
creation operators, Fock space, commutant, reflexive algebra, joint
spectral theory. }
\thanks{${}^{1\;}$partially supported by a Canadian NSERC Post-doctoral
Fellowship.}

\date{}
\begin{abstract}
Non-commutative multi-variable versions of weighted shifts arise naturally
as `weighted' left creation operators acting on Fock space. We investigate
the unital $\wot$-closed algebras they generate.
The unweighted case yields non-commutative analytic Toeplitz algebras. The
commutant can be described in terms of weighted right creation
operators when the weights satisfy a condition specific to the non-commutative
setting. We prove these algebras are reflexive when the eigenvalues for
the adjoint algebra include an open set in complex $n$-space, and provide a new
elementary
proof of reflexivity for the unweighted case. We compute eigenvalues for
the adjoint algebras in general, finding geometry not present in the
single variable
setting. Motivated by this work, we obtain general information on the
spectral theory for non-commuting $n$-tuples of operators.
\end{abstract}

\maketitle

The study of non-commutative multi-variable versions of weighted
shift operators was initiated in \cite{Kribs_simple}. These $n$-tuples arise
naturally as `weighted' left
creation operators acting on Fock space. Certain $\ca$-algebras
determined by these
{\it weighted shifts on Fock space} played a crucial role in
\cite{Kribs_simple}, and the entire class is currently under investigation
in \cite{CKM}. In this paper, we consider non-selfadjoint algebras
generated by these operators. In particular, we are interested in the weak
operator topology closed non-selfadjoint algebras they generate. There is
now an extensive body of literature for the unweighted case. The algebras
generated by the left creation operators have been established as
the appropriate {\it non-commutative analytic Toeplitz algebras}
(for instance see \cite{AP,DP1,DP2,Kribs_factor,Pop_fact,Pop_beur}).

Our motivation with this work is twofold: we wish to establish non-trivial
analogues of results obtained for standard weighted shifts. Towards this
end we are motivated by the well-known survey article \cite{Shields}. At
the same time, we wish to expose differences encountered in this new
non-commutative setting. Many of the weight conditions we obtain are
exclusive to this setting in that they reduce to trivialities in the
single variable case $n=1$.

The first section contains a review of the basic facts for these weighted
shifts, as well as an introduction to the algebras $\flambda$ we study and
the associated weight functions.
In the second section we show that, under a weight condition specific to
the non-commutative setting, a commutant theorem can be
proved which generalizes the single variable commutant theorem
\cite{Shields,SW}, as well as the unweighted $n\geq 2$ case
\cite{DP1,Pop_beur}. The
condition amounts to requiring the boundedness of particular weighted
right creation operators. This theorem leads to good
internal information on the algebras.

We investigate the reflexivity of $\flambda$ in the third section. If the
set of
eigenvalues for $\flambda^*$ includes an open set in
$\bbC^n$, then $\flambda$ is reflexive. Our proof gives a new elementary
proof
of the reflexivity of non-commutative analytic Toeplitz algebras $\fL_n$
\cite{AP,DP1}. We pose some open
questions related to the reflexivity of $\flambda$.

In the fourth section we
compute eigenvalues for $\flambda^*$, finding geometry not
present in other settings. In general there is a wealth of eigenvalues
if there are any which are non-zero. Using this analysis as motivation,
we discuss the spectral
theory for non-commuting $n$-tuples of operators in the fifth section. In
particular, we show that the natural notion of a full spectrum
\cite{Taylor1,Taylor2} does not carry over to the non-commutative several
variable case. Nonetheless, substantial information can still be
obtained from one-sided spectra.

In the final section we present examples of shifts satisfying the various
weight conditions derived throughout the paper.
Our examples include certain subclasses of the periodic shifts introduced
in \cite{Kribs_simple}. We also discover other new subclasses which help
illustrate points.

\section{Introduction}\label{S:intro}

For positive integers $n \geq 2$, let $\bbF_n^+$ be the unital free
semigroup on $n$ {\it non-commuting} letters $\{ 1, 2,
\ldots, n \}$. One way to realize $n$-variable Fock space is as the Hilbert
space $\H_n = \ell^2 (\bbF_n^+)$, where an orthonormal basis is given by vectors
corresponding to words $\{ \xi_w : w \in \bbF_n^+\}$ with the {\it vacuum
vector} $\xi_e$ corresponding to the unit or empty word $e$ in
$\bbF_n^+$. A
{\it weighted shift} on Fock space is an $n$-tuple $S = \rows$ of
operators $S_i \in \bofh$ such that there is a unitary $U: \H_n
\rightarrow \H$ and operators $T_i = U^* S_i U$ for which
there are scalars $\Lambda = \{ \lamiw \}$ with
\begin{eqnarray}
T_i \basisvector = \lamiw \xi_{iw} \qfor w\in\bbF_n^+ \qand 1\leq i \leq
n.
\end{eqnarray}

A helpful pictorial way to think of $T = \rowt$ is as a `weighted Fock
space
tree'. This was outlined in \cite{Kribs_simple}, but the basic idea is the
following: The vacuum vector $\vacuumvector$ corresponds to the vertex
lying at the top of the tree, and every other basis vector $\basisvector$
corresponds to a vertex with exactly $n$ edges leaving it downwards to the
vertices for $\xi_{iw}$, and a unique edge coming into it from above. When
we regard the edges as weighted by the scalars $\lamiw$, we can think of
this weighted Fock space tree as completely describing the actions of
$T_1, \ldots, T_n$.

As with standard weighted shifts we make some simplifying
assumptions. For the sake of brevity, we assume that the weighted shifts
$T = \rowt$ act on $\H_n$ as in $(1)$ with weights given by $\Lambda = \{
\lamiw \}$. Further, a unitary $U \in \B(\H_n)$ which is diagonal with
respect to $\{ \basisvector \}$ can be constructed for which the weighted
shift $(U^*T_1U, \ldots, U^*T_nU)$ has the non-negative weights $\{
|\lamiw|\}$. In addition,
all $T_i$ are injective precisely when each
$\lamiw \neq 0$, since there are no sinks in the weighted tree. Hence
we shall make the following assumption on $\Lambda = \{\lambda_{i,w}\}$:
\[
\mbox{{\it Assumption:     }} \lamiw > 0 \qfor w\in\fnplus \qand 1 \leq i
\leq n.
\]
We now define the algebras we wish to study.
\begin{defn}
Given a weighted shift $T = \rowt$ on $\H_n$ with weights $\Lambda = \{
\lamiw \}$, define $\fL_{\Lambda}$ to be the unital $\wot$-closed
algebra generated by the operators $\{T_1, \ldots, T_n\}$.
\end{defn}

The unweighted case $\lamiw \equiv 1$ yields the {\it left creation
operators} $L = \rowl$, of theoretical physics and free probability theory.
The algebras in this case are the so called non-commutative analytic Toeplitz algebras
$\fL_n$, which are also the $\wot$-closed algebras generated by the
left regular representation of $\fnplus$ (see
\cite{AP,DP1,DP2,Kribs_factor,Pop_fact,Pop_beur}). The operators $L = \rowl$ can be
regarded as non-commutative multi-variable versions of the unilateral
shift, hence considering weighted versions seems like a natural line of
research. The following facts were easily derived in \cite{Kribs_simple}
for $T = \rowt$:
\begin{enumerate}
\item[$(i)$] Each $T_i= L_i W_i$ where $W_i$ is the
diagonal weight operator, which is positive and injective with our above
assumption, given by $W_i \basisvector = \lamiw \basisvector$.
\item[$(ii)$] $||T_i || = \sup_w \{ \lamiw \}$, for $1 \leq i \leq n$, and
$||T|| = \sup_{i,w} \{ \lamiw \}$.
\end{enumerate}

As in \cite{CKM,CM}, it is sensible to introduce certain weight functions
when studying multi-variable weighted shifts.
\begin{defn}
Given weights $\Lambda = \{ \lamiw \}$, define the associated
{\it weight function} $W : \fnplus \times \fnplus \longrightarrow \bbR_+$
by
\[
W(u,w) = \left\{ \begin{array}{cl}
\lambda_{i_1,u} \lambda_{i_2,i_1u} \cdots \lambda_{i_k,i_{k-1}\cdots i_1u}
& \mbox{if $w=i_k \cdots i_1$} \\
1 & \mbox{if $w=e$}
\end{array}\right.
\]
\end{defn}

For words $w = i_k\cdots i_1 \in \fnplus$, we introduce the notational
convenience $T_w := T_{i_k} \cdots T_{i_1}$. Then we have
\[
T_w \xi_u = W(u,w) \xi_{wu} \qfor u,w\in\fnplus.
\]
With respect to the weighted tree structure, $W(u,w)$ is
the product of all weights picked up when one moves from the vertex
$\xi_u$ to the vertex $\xi_{wu}$. This function satisfies a cocycle
condition given by
\begin{eqnarray}
W(u,vw)  = W(wu,v) W(u,w).
\end{eqnarray}
This can be easily seen when viewing the tree structure, but we supply a
proof as well.

\begin{prop}
The weight function $W: \fnplus\times\fnplus \rightarrow \bbR_+$ satisfies
the cocycle condition $(2)$.
\end{prop}

\Prf
Let $u$, $v=i_k \cdots i_1$, and $w=j_l \cdots j_1$ be words in
$\fnplus$. Then by definition we have
\begin{eqnarray*}
W(u,vw) &=& W(u, (i_k\cdots i_1) (j_l \cdots j_1)) \\
&=& (\lambda_{j_1, u} \cdots \lambda_{j_l, j_{l-1}\cdots j_1u})
(\lambda_{i_1, wu} \cdots \lambda_{i_k, i_{k-1}\cdots i_1 wu}) \\
&=& W(u,w)W(wu,v).
\end{eqnarray*}
Finally, the formula clearly holds when $v=e$ or $w=e$.
\bx

\section{Commutant and Basic Properties of
$\fL_\Lambda$}\label{S:basicprop}

We begin by investigating the commutant structure. For $n=1$, the
operators commuting with a weighted shift are precisely the $\wot$-limits of
polynomials in the shift; in other words, the $\wot$-closed algebra
generated by the shift is its own commutant \cite{Shields,SW}.
For $n\geq 2$, the commutant of $\fL_n$ is the algebra $\fR_n$
determined by the right regular representation of $\bbF_n^+$
\cite{DP1,Pop_beur}. This is the unital $\wot$-closed algebra generated by
isometries $R_i \in \B(\H_n)$ defined by $R_i \basisvector =
\xi_{wi}$. Thus the natural generalization of these results to our setting
would require commutants determined by {\it weighted right regular
representations}. As the first lemma shows, if this is to be the case,
we are forced into specific choices for the corresponding weighted right creation
operators.

\begin{lem}\label{commlemma}
Let $T = \rowt$ be a weighted shift.
If there are operators $S = \rows$ with $S_i \in \bofhn$ which satisfy
\[
S_i \xi_w = \mu_{i,w} \xi_{wi} \qfor w\in\fnplus \qand 1\leq i \leq n,
\]
then $T_i S_j = S_j T_i$ for $1 \leq i,j \leq n$ if and only if the scalars
\begin{eqnarray}
\,\,\,\,\, \muiw = c_i W(i,w) W(e,w)^{-1} \qfor w\in\fnplus \qand 1\leq i \leq n,
\end{eqnarray}
where $c_i = \mu_{i,e}$.
\end{lem}

\Prf
For $w\in\fnplus$, a computation shows that
\[
T_iS_j \basisvector = (\mu_{j,w} \lambda_{i,wj}) \xi_{iwj}
\]
and
\[
S_jT_i \basisvector = (\lamiw \mu_{j,iw} )\xi_{iwj}.
\]
Thus the operators $T_i$ and $S_j$ commute precisely when
\[
\mu_{j,iw} = \mu_{j,w} \lambda_{i,wj} \lamiw^{-1}
\]
for all $w$ and $1\leq i,j \leq n$.
Using this formula repeatedly shows that
\[
\mu_{j,iw} = \mu_{j,e} W(j,iw) W(e,iw)^{-1}.
\]
This establishes equation $(3)$ for $|w| \geq 1$, and it clearly holds for
$w=e$.
Notice also that an analogous equation for $\lamiw$ in terms of $\muiw$
can be derived. We use this fact to generate examples in Section 6.
\bx

Given a word $w =i_1\cdots i_k$ in $\fnplus$, for operators $S_i$ as
defined in the
lemma we put $S_w:= S_{i_k} \cdots S_{i_1}$ (that the product is in
reverse order is reflective of the connection with the right regular
representation). We further define the weight function $W_\mu : \fnplus
\times \fnplus \rightarrow \bbR_+$ by
\begin{eqnarray}\label{rightweightfunction}
W_\mu (v,w) = \mu_{i_1,v} \mu_{i_2,vi_1} \cdots \mu_{i_k, vi_1\cdots
i_{k-1}},
\end{eqnarray}
Observe that
$
S_w \xi_v = W_\mu(v,w) \xi_{vw}
$
and $W_\mu (\cdot,\cdot)$ satisfies the cocycle condition
\begin{eqnarray}\label{rightcocycle}
W_\mu(u,vw)= W_\mu(u,v)W_\mu(uv,w).
\end{eqnarray}
Once again, this can be seen through the direct computation in the
following proof, or by viewing a weighted
`right' Fock space tree determined by the actions of $S_1, \ldots, S_n$.

\begin{prop}
The weight function $W_\mu : \fnplus\times\fnplus \rightarrow \bbR_+$
satisfies the cocycle condtion $(5)$.
\end{prop}

\Prf
Let $u$, $v= i_1\cdots i_k$ and $w= j_1 \cdots j_l$ be words in
$\fnplus$. Then by definition we have
\begin{eqnarray*}
W_\mu (u,vw) &=& ( \mu_{i_1,u} \cdots \mu_{i_k, ui_1\cdots i_{k-1}})
(\mu_{j_1,uv} \cdots \mu_{j_l, uvj_1\cdots j_{l-1}}) \\
&=& W_\mu(u,v) W_\mu (uv,w).
\end{eqnarray*}
Lastly, the formula clearly holds when $v=e$ or $w=e$.
\bx

We can describe the commutant of those shifts which satisfy the condition,
specific to the non-commutative multi-variable setting, which appears in
Lemma 2.1.
We note that for the rest of this section, our approach mirrors that of
Section 1 from the Davidson and Pitts paper \cite{DP1}, where the
commutant of $\fL_n$ was
computed. We wish to minimize redundancy, hence we shall leave details to
the reader when part of a proof follows the lines of
\cite{DP1}. Instead we  focus on the new aspects here.

\begin{thm}\label{commutant}
Suppose the weights $\Lambda = \{ \lamiw \}$ associated with $T = \rowt$
satisfy
\begin{eqnarray}\label{commcond}
\sup_{i,w} W(i,w) W(e,w)^{-1} < \infty.
\end{eqnarray}
For $1 \leq i \leq n$, let $S_i\in\bofhn$ be defined by
\[
S_i \basisvector = \muiw \xi_{wi} \qwhere \muiw = W(i,w) W(e,w)^{-1}.
\]
Let $\fR_\Lambda$ be the unital $\wot$-closed algebra generated by $\{S_1,
\ldots, S_n\}$.
Then the commutant of $\frightlambda$ coincides with $\flambda$.
\end{thm}

\Prf
From the lemma we have $\flambda$ contained in the commutant
$\frightlambda^\prime$. To
establish the converse inclusion
fix $A\in\frightlambda^\prime$ and set $A\xi_e = \sum_w a_w\basisvector.$
Consider the operators
\[
p_k(A) = \sum_{|w|<k} \Big(1-\frac{|w|}{k}\Big) a_w W(e,w)^{-1} T_w,
\]
which clearly belong to $\flambda$, and hence to $\frightlambda^\prime$
as well.

For $k\geq 0$, let $Q_k$ denote the projection onto $\spn \{ \basisvector
: |w| =k\}$. For $X\in\bofhn$, let
$\Phi_j(X) = \sum_{k\geq \max \{0,-j\}} Q_k X Q_{k+j}$ and put
\[
\Sigma_k (X) = \sum_{|j|<k} \Big(1 -\frac{|j|}{k}\Big) \Phi_j(X) \qfor
k \geq 1.
\]
As observed in \cite{DP1}, the generalized Cesaro sums $\Sigma_k (X)$
converge to $X$ in the strong operator topology.
But observe in this case that $\Phi_j(A) \in
\frightlambda^\prime$, since $S_iQ_k = Q_{k+1}S_i$ and hence $S_i
\Phi_j(A) = \Phi_j(A) S_i$. Thus the
sequence $\Sigma_k(A)$ belongs to
$\frightlambda^\prime$ and converges $\sot$ to $A$.

However, we also have
\[
\Sigma_k(A) \xi_e = \sum_{|w|<k} \Big( 1-\frac{|w|}{k}
\Big) a_w\basisvector = p_k(A) \xi_e.
\]
Thus for $w\in\fnplus$ it follows that
\begin{eqnarray*}
\Sigma_k(A) \basisvector &=& W_\mu(e,w)^{-1} S_w \Sigma_k(A) \xi_e \\
&=& W_\mu(e,w)^{-1} S_w p_k(A) \xi_e = p_k(A) \basisvector.
\end{eqnarray*}
Hence $A$ belongs to $\flambda$, and therefore $\flambda =
\frightlambda^\prime$ as claimed.
\bx

\begin{rem}
Notice that for $n=1$, condition (\ref{commcond}) is simply the
requirement that the
weights be bounded above, hence satisfied for all shifts.
Further, in this case the $\mu_{i,w}$ are just a constant multiple of the
weights for the shift.
It follows that Theorem~2.3 generalizes the commutant theorems discussed at the
start of the section for the unweighted $n \geq 2$ case \cite{DP1,Pop_beur}, and
for standard weighted shifts \cite{Shields,SW}.
While the class of shifts satisfying (\ref{commcond}) for $n \geq 2$ is
large, it is not all-encompassing. We provide examples in Section 6.
\end{rem}

\begin{cor}
If $\Lambda$ satisfies $(\ref{commcond})$, then $\flambda^\prime =
\frightlambda$.
\end{cor}

\Prf
This is not as trivial as the unweighted case, which follows by
symmetry. Nonetheless, it is true for the shifts satisfying
(\ref{commcond}) since the
associated $S_i$ will be bounded weighted right creation operators
commuting with the $T_i$, and hence the previous proof can be followed
along with the roles of $S_i$ and $T_j$ reversed.
\bx

There are a number of other consequences of the commutant theorem. The
proofs are simple and essentially the same as the unweighted case
\cite{DP1}, hence
we leave the details to the interested reader.

\begin{cor}\label{commcor}
For $\Lambda$ satisfying condition  $(\ref{commcond})$ we have:
\begin{enumerate}
\item[$(i)$] $\flambda$ is its own double commutant,
$\flambda^{\prime\prime} = \flambda$.
\item[$(ii)$] $\flambda$ is inverse closed.
\item[$(iii)$] The only normal elements in $\flambda$ are the scalars.
\end{enumerate}
\end{cor}

Before continuing we point out a helpful computational lemma, which shows
that elements of $\flambda$ have a generalized Fourier expansion.
We remark that for $n=1$, these Fourier expansions can be used to define
weighted $H^\infty$ and $H^2$ spaces determined by the weight sequence of
a shift \cite{Gellar1,Gellar2,Shields}. The key advantage of this approach
is that many
problems can then be phrased in terms of function theory. Some of
this theory clearly goes through for $n\geq 2$, but we do not require this
machinery here.

\begin{lem}\label{helpful}
Suppose $\Lambda$ satisfies $(\ref{commcond})$. Let $A\in\flambda$ and put
$A\xi_e =
\sum_w a_w \xi_w$. Then
\begin{eqnarray}\label{computlemma}
A \xi_v = W_\mu (e,v)^{-1} \sum_w a_w W_\mu (w,v) \xi_{wv} \qfor
v\in\fnplus.
\end{eqnarray}
\end{lem}

\Prf
From Theorem~\ref{commutant} we have
\[
A \xi_v = W_\mu (e,v)^{-1} S_v A \xi_e
= W_\mu(e,v)^{-1} \sum_w a_w W_\mu(w,v)\xi_{wv}.
\]
\upbx

For the rest of this section we assume that all $\Lambda =\{\lambda_{i,w}\}$ we
consider satisfy condition (\ref{commcond}).

\begin{prop}\label{injective}
Every non-zero $A\in\flambda$ is injective.
\end{prop}

\Prf
For non-zero $ A$ in $\flambda$ there is a word $v$ with $A\xi_v \neq 0$,
thus by
equation (\ref{computlemma}) we see $A\xi_e \neq 0$. Let $v_1$ be a word
of minimal length
for which $a_{v_1}\neq 0$. Suppose $A\xi = 0$ with $\xi = \sum_w b_w
\basisvector \neq 0$, and choose a word $v_2$ of minimal length such that
$b_{v_2}\neq 0$. Then again from Lemma~\ref{helpful}, and by the
minimality of $v_1$ and $v_2$,  we have
\begin{eqnarray*}
0 = (A\xi ,\xi_{v_1v_2}) &=& \sum_{w,u} b_w a_u W_\mu(e,w)^{-1}
W_\mu(u,w) (\xi_{uw}, \xi_{v_1v_2}) \\
&=& a_{v_1} b_{v_2}  W_\mu(e,v_2)^{-1} W_\mu (v_1,v_2) \neq 0.
\end{eqnarray*}
This contradiction shows that $A$ is injective.
\bx

Since every non-trivial idempotent has kernel we have the following.

\begin{cor}
The algebra $\flambda$ contains no non-trivial idempotents.
\end{cor}

With an extra technical condition on weights we
can prove $\flambda$ is semisimple.

\begin{thm}\label{semisimple}
If the weights $\Lambda = \{\lambda_{i,w}\}$ satisfy
\begin{eqnarray}\label{semisimpleeqn}
c_\Lambda = \inf_{v\in\fnplus} \left( W_\mu(e,v)^{-1} \liminf_{k\geq 0}
\big[ W_\mu(v,v^{k-1})\big]^{1 / k} \right) > 0,
\end{eqnarray}
then every non-zero $A \in \flambda$ has non-zero spectrum.
\end{thm}

\Prf
As before, given non-zero $A \in \flambda$ with $A\xi_e = \sum_w a_w
\xi_w$
we let $v$ be a word of minimal length such that $a_v\neq 0$. By using
equation (\ref{computlemma}) repeatedly, and the cocycle equation for
$W_\mu$, we find
\begin{eqnarray*}
A^k\xi_e &=&
\sum_{w_1, \ldots, w_k \in\fnplus} a_{w_1} \cdots a_{w_k}
\big[ W_\mu(e,w_1) W_\mu(e, w_2w_1) \cdots \\
& &  W_\mu(e, w_{k-1} \cdots
w_1) \big]^{-1}
 \big[ W_\mu(w_2,w_1) W_\mu(w_3,w_2w_1) \cdots \\
& & W_\mu(w_k,
w_{k-1}\cdots w_1) \big] \xi_{w_k \cdots w_1} \\
&=& \sum_{w_1, \ldots, w_k} a_{w_1}\cdots a_{w_k} \big[
W_\mu(e,w_1) \cdots  W_\mu(e,w_{k-1}) \big]^{-1} \\
& & W_\mu(w_k, w_{k-1}\cdots w_1) \xi_{w_k\cdots w_1} \\
&=& \left( a_v^k W_\mu(e,v)^{-k} W_\mu(v,v^{k-1}) \right)\xi_{v^k}
+ \sum_{w\neq v^k, \, |w|\geq |v|^k} b_w \basisvector.
\end{eqnarray*}
The second equality follows from the identity
\[
W_\mu(e,w_j \cdots w_1) = W_\mu( e,w_j) W_\mu(w_j, w_{j-1} \cdots w_1),
\]
for $1 \leq j < k$, which is a special case of $(5)$.
Thus condition (\ref{semisimpleeqn}) gives us an $\varepsilon > 0$ such
that the inequality
\[
||A^k||^{1 / k} \geq |(A^k\xi_e, \xi_{v^k})|^{1 / k}\geq |a_v| (c_\Lambda
- \varepsilon ) > 0,
\]
holds for all sufficiently large $k$. Therefore, it follows that $A$ has
positive spectral radius.
\bx

This result appears to be new for $n=1$.
The conclusion of the theorem is the same for $\fL_n$ \cite{DP1}, hence
together with Proposition~\ref{injective} the following consequence is
proved in the same way.
In Section 6 we present examples which satisfy both (\ref{commcond}) and
(\ref{semisimpleeqn}).

\begin{cor}
If $\Lambda$ satisfies $(\ref{commcond})$ and $(\ref{semisimpleeqn})$,
then the algebra $\flambda$ contains no
quasinilpotent elements. It follows that $\flambda$ is semisimple, and the
spectrum of every non-scalar element in $\flambda$ is connected with more
than one point.
\end{cor}

\section{Reflexivity of $\flambda$}\label{S:reflexivity}

Recall that given an operator algebra $\fA$ and a collection of subspaces
$\L$, the subspace lattice $\Lat \fA$ consists of those subspaces left
invariant by every member of $\fA$, and the algebra $\Alg \L$ consists of
all operators which leave every subspace in $\L$ invariant. Every algebra
satisfies $\fA \subseteq \Alg \Lat \fA$, and an algebra $\fA$ is {\it
reflexive} if $\fA = \Alg \Lat \fA$.

We begin with the main result of this
section. But first observe that we can regard the eigenvalues of
$\flambda^*$ as forming a subset of
$\bbC^n$, since every eigenvector $\xi$ for the adjoint algebra satisfies
equations $T_i^* \xi = \overline{\lambda}_i \xi$ for
$1\leq i \leq n$.
Given $\lambda = (\lambda_1, \ldots, \lambda_n)\in\bbC^n$ and a word
$w = i_1\cdots i_k \in \fnplus$, we shall write $w(\lambda) =
\lambda_{i_1} \cdots \lambda_{i_k} \in \bbC$.

\begin{thm}\label{reflexive}
Suppose $\Lambda$ satisfies condition $(\ref{commcond})$. If the set of
eigenvalues for
$\flambda^*$ contains an open set $\bbU$ in $\bbC^n$,  then $\flambda$
is reflexive.
\end{thm}

\Prf
Let $A \in \Alg \Lat \flambda$.
Let $\nulam$ be an eigenvector satisfying $T_i^* \nulam =
\ol{\lambda}_i\nulam$ for $1 \leq i \leq n$.
Since $\{ \nulam
\}^\perp$ belongs to $\Lat \flambda \subseteq \Lat A$, we have $A^* \nulam
=
\ol{\alpha}_\lambda \nulam$ for some scalars $\alpha_\lambda$ and all
$\ol{\lambda}= (\ol{\lambda}_1, \ldots, \ol{\lambda}_n)\in\bbU$. Let
$A\xi_e = \sum_w a_w \xi_w$, then a computation
yields
\begin{eqnarray*}
\alpha_\lambda (\xi_e,\nulam) = (A\xi_e,\nulam) &=& \sum_w a_w W(e,w)^{-1}
(\xi_e, T_w^* \nulam) \\
&=& \sum_w a_w W(e,w)^{-1} \wlam (\xi_e,\nulam).
\end{eqnarray*}
But $(\xi_e,\nulam)\neq 0$ for all $\ol{\lambda}\in\bbU$. There are a
number of
ways to see this, including an explicit formula for $\nulam$ which we
derive in the next section. Thus
$\alpha_\lambda = \sum_w a_w W(e,w)^{-1} \wlam$ for $\ol{\lambda}\in\bbU$.

Let $v \in \fnplus$, and notice that the subspace $\spn \{
\xi_{wv}: w\in\fnplus \} $
belongs to $\Lat \flambda$, hence is invariant for $A$. Thus we may write
$A\xi_v = \sum_w b_w \xi_{wv}$ for some scalars $b_w$, which gives
\begin{eqnarray*}
(A\xi_v, \nulam) &=& \sum_w b_w W(e,wv)^{-1} (T_{wv} \xi_e, \nulam) \\
&=& \sum_w b_w \left[ W(e,v) W(v,w) \right]^{-1} \wlam
v(\lambda) (\xi_e,\nulam).
\end{eqnarray*}
On the other hand, we have
\begin{eqnarray*}
(A\xi_v, \nulam) &=& W(e,v)^{-1} (\xi_e, T_v^* A^* \nulam) \\
&=& W(e,v)^{-1} \alpha_\lambda v(\lambda) (\xi_e,\nulam).
\end{eqnarray*}
It follows that $\alpha_\lambda = \sum_w b_w W(v,w)^{-1} \wlam$ for
$\ol{\lambda}\in\bbU$.
We shall see in Theorem 4.3 and its corollaries, that the hereditary
nature of the eigenvalue set for $\flambda^*$ allows us to assume that the
origin in $\bbC^n$ belongs to $\bbU$. Thus, as $\bbU$ is an open set
containing the origin, and hence contains a polydisc about the origin, the
theory of analytic functions of several variables \cite{Rudin} gives us
$b_w W(v,w)^{-1} = a_w W(e,w)^{-1}$ for $w\in\fnplus$. In particular,
\[
A\xi_v = \sum_w a_w W(v,w)W(e,w)^{-1} \xi_{wv}.
\]

Recall from Theorem~\ref{commutant} that $\flambda =
\frightlambda^\prime$, and $\frightlambda$ is generated by $S_1, \ldots
,S_n$ (and the unit). Using the definition of $S_i$ we get:
\[
(AS_i) \xi_v =  \mu_{i,v} A\xi_{vi} =  \sum_w a_w \left( \mu_{i,v}
W(vi,w)W(e,w)^{-1} \right) \xi_{wvi}
\]
and
\[
(S_iA) \xi_v = \sum_w a_w \left( W(v,w) W(e,w)^{-1}
\mu_{i,wv}\right) \xi_{wvi}.
\]
But from the definition of $\mu_{i,w}$, and the cocycle equation $(2)$ for
$W(\cdot,\cdot)$, we see that
\begin{eqnarray*}
\mu_{i,v} W(vi,w) &=& W(i,v) W(e,v)^{-1} W(vi,w) \\
&=& W(e,v)^{-1} W(i,wv) \\
&=& W(e,wv)^{-1} W(v,w) W(i,wv) \\
&=& \mu_{i,wv} W(v,w)
\end{eqnarray*}
Heuristically, both sides of this equation consist of the product of
weights from $i$ to $wvi$ in the weighted Fock space tree, divided by the
product of weights from $e$ to $v$.
Therefore, $AS_i = S_iA$ for $1 \leq i \leq n$, and $A$ belongs to
$\frightlambda^\prime = \flambda$, as required.
\bx

\begin{rem}
This gives a new elementary proof of the reflexivity of $\fL_n$ for
$n\geq 2$. Indeed, in \cite{AP} reflexivity is proved using
non-commutative factorization results for $\fL_n$, and the authors point
out that their proof does not carry over to the commutative case. Further,
in \cite{DP1} the stronger notion of hyper-reflexivity is proved for
$\fL_n$, hence subsuming reflexivity. But as the authors point out, it is
not clear how just reflexivity follows from their proof without using the
stronger notion,
whose proof uses some deep facts from von Neumann algebra theory. In the
event, our proof simply relies on several variable analytic function
theory in a way which carries over to the commutative case $n=1$.
We also mention that the paper \cite{KP} on {\it free semigroupoid
algebras} contains a different elementary proof of
reflexivity for $\fL_n$.

The hypothesis of the theorem can be weakened to simply
require
an eigenvalue $\lambda = (\lambda_1, \ldots, \lambda_n)$ for $\flambda^*$
with each $\lambda_i \neq 0$. For if this is the case,
Corollary~\ref{openset} shows there is also an open set of eigenvalues
about the origin.

For $n=1$, every weighted shift $T$ for which $T^*$ has
a non-zero eigenvalue is reflexive. In other words, the algebra  generated
by $T$ is reflexive \cite{Shields}. However, when there is a non-zero
eigenvalue
for $T^*$, the spectral theory for unilateral weighted shifts implies
there is an
entire disc of such eigenvalues. Thus Theorem~\ref{reflexive} is a
generalization of this result. By using
Corollary~\ref{bdedbelowcor}, we obtain
examples which fulfil the hypothesis of the theorem in Section 6.
\end{rem}

At the present time we have no examples of non-reflexive algebras $\flambda$ for $n\geq 2$.
There are numerous examples in the single variable case; for instance, every unicellular
weighted shift operator \cite{Harrison,Nik,Shields} generates a non-reflexive algebra. This
motivates the following problem.

\begin{prob}
Is there an analogue of unicellular weighted shifts in the non-commutative multivariable
setting?
\end{prob}

For $v\in\fnplus$, let $E_v$ be the subspace $E_v = \spn \{ \xi_{wv}
: w\in\fnplus \}$. If $T
= \rowt$ is a weighted shift, then each
$E_v$ is clearly invariant for each $T_i$, and thus belongs to $\Lat
\flambda$.
A natural guess for a generalization of unicellular weighted shifts here would be $\Lat \flambda$
equal to the subspace lattice $\L$ generated by $\{ E_v : v\in\fnplus\}$. But there are simple
examples which show this cannot happen for $n\geq 2$. Indeed, the subspace $\M = \spn \{
\xi_e, \frac{1}{\sqrt{2}}(\xi_1 + \xi_2)\}$ clearly belongs to $\Lat \flambda^*$ for all
$\Lambda$, and hence $\M^\perp$ is a subspace outside of $\L$ which belongs to every $\Lat
\flambda$.

More generally,
it would also be interesting to characterize when $\flambda$ is
reflexive. This appears to be an open problem for
$n=1 $ still.

\begin{prob}
Determine which algebras $\flambda$ are reflexive.
\end{prob}

\section{Eigenvalues for $\flambda^*$}\label{S:eigenvalues}

It is clear that unilateral weighted shifts  have no non-zero eigenvalues,
and it is also obvious that the
$T_i$ from a shift $T = \rowt$ have no non-zero eigenvalues. However, the
adjoint of a unilateral weighted shift has
a wealth of non-zero eigenvalues if it has any which are non-zero. We
mention that for $n=1$ the point spectrum for
the adjoint of a weighted shift is either an open disc, possibly with some
boundary points included, or just the origin \cite{Kelley,Shields}. For
$n\geq 2$, the eigenvalues for $\fL_n^*$ form the open complex unit $n$-ball
$\bbB_n = \{ \lambda\in\bbC^n : ||\lambda||_2 < 1 \}$ \cite{AP,DP1}.
In the eigenvalue analysis of this section we discover geometry not
present in either of these settings.
Let us begin the discussion by observing some basic facts.

\begin{prop}
If $T = \rowt$ is a weighted shift, then $T_i^*\xi_e = 0$ for $1 \leq i
\leq n$, and hence $\vec{0}\in\bbC^n$ is always an eigenvalue for
$\flambda^*$. In fact more is true,
\[
\bigcap_{i=1}^n \ker T_i^* = \spn \{ \xi_e \}.
\]
\end{prop}

\Prf
It is clear that $T_i^* = W_i L_i^*$ annihilates the vacuum vector. If
$\xi$ satisfies $T_i^* \xi = W_i L_i^* \xi = 0$, then $L_i^* \xi = 0$
since $W_i$ is injective. But the projection $P_e =\xi_e \xi_e^*$ onto the
vacuum space
is given by the equation $P_e = I - \sum_{i=1}^n L_i L_i^*$. Hence
$P_e \xi = \xi$ and $\xi$ belongs to $\spn \{ \xi_e \}$.
\bx

More generally, if $\lambda = (\lambda_1, \ldots, \lambda_n)\in\bbC^n$ and
some unit vector $\xi$ in $\H_n$ satisfy $T_i^* \xi = \ol{\lambda}_i \xi$,
then
\[
||\lambda||_2^2 = \sum_{i=1}^n |\lambda_i|^2 = \sum_{i=1}^n
(T_iT_i^* \xi,\xi)
= (TT^* \xi,\xi)  \leq ||T||^2 = \sup_{i,w} \lamiw^2.
\]
Thus all eigenvalues for $\flambda^*$ are contained in the corresponding
ball in $\bbC^n$ determined by the supremum of $\Lambda = \{ \lamiw
\}$. In general this
estimate is not a good one though.

We now determine the form of all joint eigenvectors. Notice we make no
extra assumptions on weights in the next three results.

\begin{lem}\label{evectorlemma}
Let $T = \rowt$ be a weighted shift with weights $\Lambda$. If a vector
$\xi= \sum_w \ol{a}_w \xi_w \in \H_n$ satisfies $T_i^* \xi =
\ol{\lambda}_i \xi$ for $1 \leq i \leq n$, then
\begin{eqnarray}\label{evectoreqn}
a_w = a_e \wlam W(e,w)^{-1} \qfor w\in\fnplus.
\end{eqnarray}
Conversely, if $\xi$ belongs to $\H_n$ with coefficients
$\ol{(\xi,\basisvector)}$ satisfying $(\ref{evectoreqn})$, normalized so
that $a_e = (\xi,\xi_e) = 1$, then $T_i^* \xi
= \ol{\lambda}_i\xi$ for $1 \leq i \leq n$. It follows that every
eigenspace for $\flambda^*$ is
one-dimensional.
\end{lem}

\Prf
If $\xi = \sum_w \ol{a}_w \xi_w$ is in $\H_n$ with
$T_i^* \xi = \ol{\lambda}_i \xi$, then
\[
\ol{\lambda_i a_e} = (\ol{\lambda}_i \xi, \xi_e) = (T_i^* \xi, \xi_e)
= \lambda_{i,e}(\xi,\xi_i) = \ol{a}_i \lambda_{i,e},
\]
so that $a_i = a_e \lambda_i \lambda_{i,e}^{-1}$. Another computation
gives the next level:
\[
\ol{\lambda_i a_j} = (\ol{\lambda}_i \xi, \xi_j) = (T_i^*\xi, \xi_j) =
\lambda_{i,j} (\xi,\xi_{ij}) =\ol{a}_{ij} \lambda_{i,j},
\]
hence $a_{ij} = a_j \lambda_i \lambda_{i,j}^{-1}
= a_e (\lambda_i \lambda_j)(\lambda_{i,j} \lambda_{j,e})^{-1}.$ We can now
obtain equation (\ref{evectoreqn}) by induction.

On the other hand, if $\xi = \sum_w \ol{\wlam} W(e,w)^{-1} \xi_w$ belongs
to $\H_n$, then since $T_i^* = W_i L_i^*$ we have
\begin{eqnarray*}
T_i^* \xi = W_i L_i^* \xi &=& \sum_{w=iu; \, u\in\fnplus} \ol{\wlam}
W(e,w)^{-1} W_i L_i^* \xi_w \\
&=& \sum_u \ol{(iu)(\lambda)} W(e,iu)^{-1} \lambda_{i,u} \xi_u  =
\ol{\lambda}_i \xi.
\end{eqnarray*}
It also follows that for every eigenvalue $\lambda = (\lambda_1 ,\ldots,
\lambda_n)$ of $\flambda^*$, the eigenspace $E_\lambda = \spn \{ \xi :
T_i^* \xi = \lambda_i \xi \qfor 1\leq i \leq n \}$ is one-dimensional.
\bx

As an immediate consequence we get a
tight characterization of the eigenvalues for $\flambda^*$ which
generalizes the known result for $n=1$ \cite{Gellar2}.

\begin{thm}\label{evaluethm}
The eigenvalues for $\flambda^*$ consist of all $\lambda\in\bbC^n$ for
which
\begin{eqnarray}\label{evalueseries}
\sum_{w\in\fnplus} |\wlam|^2 W(e,w)^{-2} < \infty.
\end{eqnarray}
\end{thm}

In particular, it is evident that the set of eigenvalues has the following
hereditary property.

\begin{cor}
If $\mu = (\mu_1, \ldots, \mu_n)$ is an eigenvalue for $\flambda^*$, then
so is every element of the set
\[
\bbD = \{ (\lambda_1, \ldots, \lambda_n) : |\lambda_i| \leq |\mu_i| \qfor
1 \leq i \leq n\}.
\]
\end{cor}

\Prf
If $\lambda$ belongs to $\bbD$, then the series
$\sum_w |w(\lambda)|^2 W(e,w)^{-2}$ converges since its partial sums are
bounded above by the corresponding series for $\mu$, and
the result follows from the theorem.
\bx

As a special case of this result, together with Theorem~\ref{reflexive},
we discover a method for finding reflexive algebras.

\begin{cor}\label{openset}
If $\Lambda$ satisfies condition $(\ref{commcond})$, and if $\flambda^*$
has
an eigenvalue $\mu= (\mu_1, \ldots, \mu_n)$ with each
$\mu_i \neq 0$, then the eigenvalues for $\flambda^*$ include a polydisc
in $\bbC^n$, and hence $\flambda$ is reflexive.
\end{cor}

Let us be more concrete. Specifically, when the weights $\Lambda=
\{\lambda_{i,w}\}$ for a shift are bounded away from zero
we identify ellipses as eigenvalue sets.

\begin{cor}\label{polyellipsecor}
If  $\Lambda= \{\lambda_{i,w}\}$ is bounded away from zero, then
the eigenvalues for $\flambda^*$ include the set
\[
\bbE = \Big\{ (\lambda_1, \ldots , \lambda_n)\in\bbC^n \,\,\Big| \,\,
\frac{|\lambda_1|^2}{c_1^2} + \ldots + \frac{|\lambda_n|^2}{c_n^2} < 1
\Big\},
\]
where $c_i = \inf_{w\in\fnplus} \lamiw > 0$ for $1\leq i \leq n$.
\end{cor}

\Prf
Let $c=(c_1,\ldots, c_n)$.
Then for $w= i_k \cdots i_1\in\fnplus$ we have
\[
w(c) =c_{i_k} \cdots c_{i_1}  \leq \lambda_{i_k,i_{k-1}\cdots i_1}
\cdots \lambda_{i_2,i_1} \lambda_{i_1,e} = W(e,w).
\]
Let $\lambda = (\lambda_1, \ldots, \lambda_n)$ belong to $\bbE$
and put $\lambda c^{-1} = (\lambda_1 c_1^{-1}, \ldots, \lambda_n
c_n^{-1})$.
Then $||\lambda c^{-1} ||_2 < 1$, and hence
\begin{eqnarray*}
\sum_{w\in\fnplus} |\wlam|^2 W(e,w)^{-2} &=& \sum_w |w(\lambda c^{-1})|^2
w(c)^2
W(e,w)^{-2} \\
&\leq& \sum_w |w(\lambda c^{-1})|^2 \\
&=& \sum_{k\geq 0} \left( \sumin |\lambda_i|^2 c_i^{-2} \right)^k \\
&=& (1 - ||\lambda c^{-1} ||_2^2)^{-1} < \infty.
\end{eqnarray*}
Thus by Theorem~\ref{evaluethm}, every $\lambda = (\lambda_1, \ldots,
\lambda_n)$ in $\bbE$ is an eigenvalue for $\flambda^*$.
\bx

Therefore, with Theorem~\ref{reflexive} this result gives us a large
subclass of reflexive algebras.
We give examples of such shifts in Section 6.

\begin{cor}\label{bdedbelowcor}
If $\Lambda= \{\lambda_{i,w}\}$ is bounded away from zero and satisfies
$(\ref{commcond})$, then $\flambda$ is reflexive.
\end{cor}

\section{Joint Spectral Theory}\label{S:spectral}

Our analysis in the previous section gives us motivation for discussing
the joint spectral theory for non-commuting
$n$-tuples of operators. We mention the seminal work of Taylor
\cite{Taylor1,Taylor2} on multi-variable spectral theory.
As suggested at the end of the paper \cite{Taylor1}, a natural definition
for a (non-commuting) $n$-tuple of operators $A= \rowa$ on $\H$ to be
non-singular would be to require the existence of another $n$-tuple
$B = (B_1, \ldots ,B_n)$ for which
\begin{eqnarray}\label{rightinv}
A_1B_1 + \ldots + A_nB_n =I
\end{eqnarray}
and
\begin{eqnarray}\label{leftinv}
B_i A_j =\delta_{ij} I \qfor 1 \leq i,j \leq n.
\end{eqnarray}
This is a natural notion since it says the map $A:\H^{(n)} \rightarrow \H$
given by $A [x_1 \cdots x_n]^t = \sum_{i=1}^n A_i x_i$ is
invertible. Indeed, equations (\ref{rightinv}) and (\ref{leftinv}) can be
written as $AB^t = I_\H$ and $B^t A = I_{\H^{(n)}}$.

Hence we shall use the
following nomenclature: the {\it joint right spectrum} $\sigma_r (A)$ of
$A= \rowa$ consists of $\lambda \in \bbC^n$ for which no solution
$B = (B_1, \ldots ,B_n)$ exists to the equation
\[
\sum_{i=1}^n (A_i - \lambda_i I)B_i =I.
\]
Whereas the {\it joint left spectrum} $\sigma_l (A)$ consists of
$\lambda\in\bbC^n$ for which no solution $B= (B_1, \ldots, B_n)$ exists to
the equations
\[
B_i (A_j - \lambda_jI) = \delta_{ij}I \qfor 1 \leq i,j \leq n.
\]
The {\it full spectrum} $\sigma(A)$ is the union of $\sigma_r(A)$ and
$\sigma_l(A)$.

Our goal in this section is to show that, while the right spectrum can
yield good information, the above notions of left and
full spectra are inappropriate for the non-commutative multi-variable
setting. We accomplish this by focusing on the unweighted case. Let
us begin by computing the right spectrum of the left creation operators.

\begin{thm}\label{rightspectrum}
The right spectrum of $L = \rowl$ is $\sigma_r(L) = \ol{\bbB}_n$.
\end{thm}

\Prf
From Corollary~\ref{polyellipsecor}, it follows that every
$\lambda\in\bbB_n$
determines an eigenvector $\nulam = \sum_w \ol{\wlam} \basisvector$ for
$\fL_n^*$ (The eigenvalues for the unweighted case $\fL_n^*$ were
initially worked out in \cite{AP,DP1}.). Thus for $\lambda\in\bbB_n$ there
can be no solution $B=
(B_1, \ldots, B_n)$ to
$\sum_{i=1}^n B_i^*(L_i^* - \ol{\lambda}_iI) =I$, so that $\bbB_n
\subseteq \sigma_r(L)$, and a standard approximation argument shows that
$\ol{\bbB}_n \subseteq \sigma_r(L)$.

On the other hand, for $\lambda\in\bbB_n$ the operator $I- \sum_{i=1}^n
\ol{\lambda}_iL_i$ is invertible since
\[
\big|\big| \sum_i \ol{\lambda}_i L_i \big| \big|^2 = \sum_i |\lambda_i|^2
=
||\lambda||^2 < 1.
\]
A computation shows that $(I-\sum_{i=1}^n \ol{\lambda}_i L_i)^{-1} =
\sum_w \ol{\wlam} L_w$ as linear transformations, hence the latter
operator is bounded and the
sum is a $\wot$-limit. Thus if
$\lambda \in \bbC^n$ with $||\lambda|| >1$, and defining
$\lambda^{-1} = \frac{1}{||\lambda||^2} (\ol{\lambda}_1, \ldots,
\ol{\lambda}_n)$, we have $\sum_w \ol{w(\lambda^{-1})}L_w$ belonging to
the inverse closed algebra $\fL_n$. For $1 \leq i \leq n$, define
\[
B_i = \frac{\ol{\lambda}_i}{||\lambda||^2} \sum_w \ol{w(\lambda^{-1})}
L_w.
\]
It follows that
\[
\sumin (\lambda_i I - L_i) B_i = \sum_w \ol{w(\lambda^{-1})} L_w -
\sumin \Big( \sum_{u=iw, w\in\fnplus} \ol{u(\lambda^{-1})}L_u\Big)  = I,
\]
and hence $\lambda\in\sigma_r(L)^c$. Therefore, $\sigma_r(L) =
\ol{\bbB}_n$ as claimed.
\bx

This theorem generalizes the well-known result for $n=1$, that the right
spectrum of the unilateral weighted shift operator is the closed unit
disc. Of course, in the single variable case we know this is the entire
spectrum.
The rigidity in the definition of left invertibility for $n\geq 2$,
results in $L = \rowl$ having a very large left spectrum.

\begin{prop}\label{leftspectrum}
Let $\bbD_n = \{ \lambda = (\lambda_1, \ldots, \lambda_n) \in\bbC^n
: |\lambda_i|<1\}$ be the unit polydisc in $\bbC^n$.
The left spectrum of $L = \rowl$ includes the complement of this set,
\[
\bbC^n \setminus \bbD_n \subseteq \sigma_l(L).
\]
\end{prop}

\Prf
For the sake of brevity, we focus on the $n=2$ case. We are required to
show there are no solutions $A= (A_1, A_2)$ to
\begin{eqnarray}\label{creationeqn}
A_i (L_j - \lambda_j I) = \delta_{ij} I \qfor 1 \leq i,j \leq 2,
\end{eqnarray}
when $|\lambda_1| \geq 1$ or $|\lambda_2|\geq 1$.
Suppose there was such a solution. Then $A_2 L_1 = \lambda_1 A_2$,
and if $|\lambda_1|>1$ we would have
\[
A_2 \xi_{1^k} = \lambda_1 A_2 \xi_{1^{k-1}} = \ldots = \lambda_1^k A_2
\xi_e \qfor
k\geq 1
\]
implying the unboundedness of $A_2$ unless $A_2 \xi_e=0$, which is
addressed below.

If instead $|\lambda_1|=1$, then $A_1
L_1
=\lambda_1 A_1 + I $, and for $k \geq 1$ we would have
\begin{eqnarray*}
A_1 \xi_{1^k} = A_1 L_1 \xi_{1^{k-1}} &=& \lambda_1 A_1 \xi_{1^{k-1}} +
\xi_{1^{k-1}} \\
&\vdots&   \\
&=& \lambda_1^k A_1 \xi_e + \lambda_1^{k-1} \xi_e + \lambda_1^{k-2}\xi_1 +
\ldots + \xi_{1^{k-1}}.
\end{eqnarray*}
Hence $\lim_{k \rightarrow \infty} ||A_1 \xi_{1^k} || = \infty$, implying
the unboundedness of $A_1$.
Similarly, $|\lambda_2| =1$ implies that $A_2$ is unbounded.
Finally, if $|\lambda_2|>1$ and $A_2\xi_e=0$, then the equation
$A_2L_2 = \lambda_2 A_2 + I$ can be applied in the same manner to obtain
$\lim_{k \rightarrow \infty} ||A_2 \xi_{2^k}|| = \infty$. Thus, in all
cases $A_1$ or $A_2$ would have to be unbounded.
\bx

\begin{note}
Notice there can be solutions to (\ref{creationeqn}) for certain values of
$\lambda = (\lambda_1, \lambda_2)$. For instance, the solutions $A = (A_1,
A_2)$ for $\lambda = \vec{0}$ are rank one perturbations of $(L_1^*,L_2^*)$
of the form
\[
A_i = L_i^* + \eta_i \xi_e^* \qfor \eta_i \in \H_2.
\]
Indeed, these are solutions to
(\ref{creationeqn}) for $\lambda = \vec{0}$ since $(L_1^*,L_2^*)$ is and
\[
(\eta_i \xi_e^*)L_j = \eta_i (L_j^*\xi_e)^* = 0 \qfor 1\leq i,j \leq 2.
\]
Conversely, if $A=(A_1,A_2)$ is a solution for $\lambda = \vec{0}$, then
\begin{eqnarray*}
A_1 &=& A_1 L_1L_1^* + A_1L_2L_2^* + A_1 P_e \\
&=& L_1^* + (A_1 \xi_e) \xi_e^*,
\end{eqnarray*}
and the analogue is true for $A_2$.

It should be possible to say more about the left spectrum of $L$. In fact,
we believe $\sigma_l (L)$ includes $\bbC^n \setminus \ol{\bbB}_n$, in
other
words that there are {\it no} solutions to (\ref{creationeqn}) for
$||\lambda_2||>1$. This would imply that the full spectrum of $L$ is all
of $\bbC^n$. Nonetheless, we can say the following.
\end{note}

\begin{cor}
The full spectrum of $L = \rowl$  includes the set
\[
\ol{\bbB}_n \cup (\bbC^n \setminus \bbD_n) \subseteq \sigma(L).
\]
\end{cor}

\begin{rem}
Thus, towards a joint spectral theory for non-commuting $n$-tuples of
operators, it seems unreasonable to consider the definition for the full
spectrum discussed at the start of the section. Especially considering the
fact that $L= \rowl$ is fundamental
to the theory of non-commuting $n$-tuples of operators, playing the role
of the
unilateral shift in this setting. However, we have seen that the right
spectrum can yield good information for an $n$-tuple.
We mention another possibility for the left spectrum might
be to use the equation $BA^t = \sumin B_iA_i = I$ to define left
invertibility. For instance,
it appears that with this definition it may be possible to
generalize the result from \cite{Ridge} for the left spectra of weighted
shift operators to this non-commutative setting.
\end{rem}

\section{Examples}\label{S:examples}

In this section we describe examples of shifts satisfying the various
weight conditions discovered in the paper. We begin by presenting a simple
subclass satisfying (\ref{commcond}).

\begin{eg}\label{commutanteg2}
The set of shifts satisfying (\ref{commcond}) is clearly closed under
changing
finitely many weights. For instance, the shifts
\[
T = (c_1L_1, \ldots, c_nL_n), \qfor c_i > 0
\]
are easily seen to satisfy (\ref{commcond}). In fact since the $T_i$
are just multiples of $L_i$, the algebra they generate is $\fL_n$. Let
$\lamiw
> 0$ for $1\leq i \leq n$ and $|w| \leq k$, and put $\lamiw = c_i$ for
$|w|>k$. Then the shift $\td{T}$ determined by $\Lambda = \{ \lamiw \}$
will satisfy (\ref{commcond}). For, if $|w|>k$ with $w= i_s \cdots i_1$,
then
\[
\lambda_{i_j, i_{j-1}\cdots i_1} = c_{i_j} = \lambda_{i_j, (i_{j-1}\cdots
i_1) i} \qfor k< j \leq s \qand 1\leq i \leq n.
\]
Thus the weight products $W(e,w)$ and $W(i,w)$ will be determined by the
original $n$-tuple when we move far enough down the tree.

Observe that the
operators which determine the $n$-tuple $\td{T}$ will be particular
finite rank perturbations of the $c_i L_i$ in this case.
Further, the weights are bounded away from zero in this class, hence by
Corollary~\ref{bdedbelowcor} the associated algebras $\flambda$ are
reflexive.
\end{eg}

We next show that condition (\ref{commcond}) is not satisfied by all shifts.

\begin{eg}\label{commutanteg1}
For simplicity consider $n = 2$. Let $T= (T_1, T_2)$ be a weighted shift
which satisfies
\[
\lambda_{1,1^k} = m^{-1} \qand \lambda_{1,1^k2} = m^{-1 / 2} \qfor k \geq
0.
\]
Then $W(e,1^k) = m^{-k}$ and $W(2,1^k) = m^{-k / 2}$, so that
\[
\sup_{i,w} W(i,w) W(e,w)^{-1} \geq \sup_{k\geq 0} W(2,1^k)W(e,1^k)^{-1} =
\lim_{k\rightarrow\infty} m^{k / 2}.
\]
Thus if $m > 1$, equation (\ref{commcond}) is not satisfied.

We can also use this example to find shifts which do satisfy
(\ref{commcond}). Indeed, consider the above weights with $m\leq 1$, and
define all
other weights $\lamiw \equiv c$ for some constant $c \leq \sqrt{m}$.
If $w =1^k$, we have $W(1,w)=m^{-k}$, so that
\[
W(i,1^k)W(e,1^k)^{-1} \leq \max \{ 1, m^{k / 2} \} \leq 1
\]
for $i=1,2$, by the previous paragraph.
Here are the other cases for a word $w$ of length $k$:
\[
W(e,w) = \left\{ \begin{array}{cl}
m^{-s}c^{k-s} & \mbox{if $w= u 2 \, 1^s$, with $s>0$} \\
m^{-s / 2} c^{k-s} & \mbox{if $w=1^s2$ or $w= u2\, 1^s 2$, with $s\geq 0$}
\end{array}\right.
\]
Whereas, for $i=1,2$ we have
\[
W(i,w) = \left\{ \begin{array}{cl}
m^{-s}c^{k-s} & \mbox{if $w = u 2 \, 1^s$, with $s>0$ and $i=1$} \\
m^{-s / 2}c^{k-s} & \mbox{if $w = u2\, 1^s$, with $s> 0$ and $i=2$} \\
c^k & \mbox{if $w= 1^s 2$ or $w= u2\, 1^s 2$, with $s\geq 0$}
\end{array}\right.
\]
This exhausts all cases, hence it follows that the supremum in
(\ref{commcond}) is equal to $1$.
Further, the weights in this class of examples are bounded away from zero,
hence by
Corollary~\ref{bdedbelowcor} the associated algebras $\flambda$ are
reflexive.
\end{eg}

The periodic weighted shifts introduced in \cite{Kribs_simple} provide a
useful subclass for the purposes here.

\begin{eg}\label{periodic}
For $k\geq 1$, a weighted shift $T= \rowt$ is of {\it period} $k$ if for
all $w\in\fnplus$ we have
\[
T_i \xi_w = \lambda_{i,u} \xi_{iw},
\]
where $w=uv$ is the unique decomposition of $w$ with $0\leq |u| < k$ and
$|v| \equiv 0 \, ({\rm mod}\,\,k)$. The {\it remainder} scalars $\{
\lambda_{i,u} : 0 \leq |u| < k\}$ completely determine the shift. It is
most satisfying to think of this notion of periodicity in terms of the
weighted Fock space tree, which is generated by the finite weighted tree
top with vertices $\{ \xi_w : |w|\leq k\}$ and weighted edges given by the
remainder scalars. In \cite{Kribs_simple}, this notion of periodicity was
used to find non-commutative generalizations of the Bunce-Deddens
$\ca$-algebras.

Notice the 1-periodic shifts are simply of the form $(c_1L_1, \ldots
,c_nL_n)$, but for $k\geq 2$ the $k$-periodic shifts give us a large
non-trivial class to work with.
For the sake of brevity, we focus on the 2-periodic shifts
$T=(T_1,T_2)$; that is, we consider $k=2$ and $n=2$. Each such
2-tuple will be determined by six scalars, which we denote by
\[
\lambda_{1,e} = a,\,\, \lambda_{2,e}=b,\,\, \lambda_{1,1}=c,\,\,
\lambda_{2,1}=d,\,\,
\lambda_{1,2}= e,\,\, \lambda_{2,2} =f.
\]

We first observe that periodic shifts can fail to satisfy
(\ref{commcond}). For
$k\geq 1$, consider the words $w_k = (21)^k \in \bbF_2^+$. Then by
2-periodicity we have
\[
W(e,w_k) = (ad)^k \qand W(2,w_k) = (eb)^k.
\]
Consequently, if $eb > ad >0$, then
\[
\sup_{i,w} W(i,w)W(e,w)^{-1} \geq \sup_{k\geq 0} W(2,w_k)W(e,w_k)^{-1} =
\infty,
\]
showing that (\ref{commcond}) fails in this case. Nevertheless, there are
many periodic shifts which do satisfy (\ref{commcond}).

The set $\Lambda$ associated with
every periodic shift with non-zero weights is plainly bounded away from zero.
Hence by
Corollary~\ref{bdedbelowcor}, the algebra $\flambda$ generated by a
periodic shift which satisfies (\ref{commcond}) is reflexive. As a simple
subclass of examples, consider the case when
$a=b>0$ and $c=d=e=f>0$. In this case we compute
\[
W(e,w) = \left\{ \begin{array}{cl}
a^kc^k & \mbox{if $|w| = 2k$} \\
a^{k+1} c^k & \mbox{if $|w| = 2k+1$}
\end{array}\right.
\]
Whereas, for $i=1,2$ we have
\[
W(i,w) = \left\{ \begin{array}{cl}
c^ka^k & \mbox{if $|w|=2k$} \\
c^{k+1} a^k & \mbox{if $|w|=2k+1$}
\end{array}\right.
\]
Thus (\ref{commcond}) is clearly satisfied for the entire subclass,
with $ca^{-1}$ providing an upper bound, and
hence all of the associated algebras $\flambda$ are reflexive.
\end{eg}

The next observation follows directly from the analysis in Section 2,
and will allow us to generate non-trivial examples
satisfying (\ref{semisimpleeqn}).

\begin{prop}\label{rightcreation}
Let $S = \rows$ be operators in $\B(\H_n)$ defined by $S_i \xi_w = \muiw
\xi_{wi}$ for scalars
$\muiw > 0$, normalized so that each $\mu_{i,e}=1$. Let
$\td{W} : \fnplus \times \fnplus \rightarrow \bbR_+$ be the weight
function defined by $\td{W} (v,e)=1$ for all $v\in\fnplus$ and
\[
\td{W} (v,w) = \mu_{i_1, v} \mu_{i_2,vi_1} \cdots \mu_{i_k,vi_1\cdots
i_{k-1}} \qif w = i_1 \cdots i_k.
\]
Suppose that
\begin{eqnarray}\label{tilde}
\sup_{i,w} \td{W} (e,w)^{-1} \td{W} (i,w) < \infty.
\end{eqnarray}
Let $T_i \in \bofhn$ be defined by $T_i \xi_w = \lamiw \xi_{iw}$ where
\[
\lamiw = \td{W}(i,w) \td{W}(e,w)^{-1}.
\]
Then $T = \rowt$ satisfies
$(\ref{commcond})$, and $S = \rows$ are the weighted right creation
operators obtained for $T$ as in Theorem~\ref{commutant}. In particular,
the functions $\td{W} = W_\mu$ are the same.
\end{prop}

\Prf
By following the proof of Lemma~\ref{commlemma}, we see that the
$S_i$ and $T_j$ commute, and we have the corresponding formulas for
$\lamiw$ and $\muiw$ in terms of the other. The $n$-tuple $T = \rowt$
satisfies
(\ref{commcond}) precisely because the $\muiw$ are uniformly bounded. Thus
the $S_i$ really are the weighted right creation operators obtained in
Theorem~\ref{commutant} which generate the commutant, and hence $\td{W} =
W_\mu$.
\bx

\begin{eg}\label{semisimpleeg}
Let $S= \rows$ be weighted right creation operators with non-zero weights
$\muiw$ for which
there is a $k \geq 0$ with $\muiw \equiv 1$ for $|w| > k$ (There is
complete freedom on weight choices for $|w| \leq k$.). Since only finitely
many weights are different than 1, it is easy to see that $W_\mu = \td{W}$
satisfies (\ref{semisimpleeqn}) and (\ref{tilde}).

Thus by Proposition~\ref{rightcreation}, the corresponding weighted left
creation operators $T =\rowt$ satisfy (\ref{commcond}), and hence the
associated
algebras $\flambda$ are semisimple. We note that, while $\flambda$ has
this structure, the weights
$\Lambda = \{ \lamiw \}$ are not easily described. Indeed, an examination
of the formula $\lamiw = W_\mu(i,w) W_\mu (e,w)^{-1}$ shows
that typically
these scalars will not satisfy a finiteness condition analogous to the one
which defines $\muiw$.
In other words, the $S_i$ will be finite rank perturbations of the $R_i$,
but the $T_i$ will not in general be finite rank perturbations of the
$L_i$.
\end{eg}


{\noindent}{\bf Acknowledgements.}
The author is grateful to Raul Curto, Paul Muhly, and
Stephen Power for enlightening discussions. Thanks also to the
Department of Mathematics and Statistics at Lancaster University for kind
hospitality during the preparation of this article.



\begin{tabbing}
{\it E-mail address}:xx\= \kill \noindent {\footnotesize\it
Mailing Address}:
\>{\footnotesize\sc Department of Mathematics and Statistics}\\
\>{\footnotesize\sc University of Guelph}\\
\>{\footnotesize\sc Guelph, ON}\\
\>{\footnotesize\sc CANADA \quad N1G 2W1}\\
\\
{\footnotesize\it E-mail address}: \>{\footnotesize\sf
dkribs@uoguelph.ca}
\end{tabbing}

\end{document}